\documentclass[11pt]{amsart}
\usepackage[english]{babel}
\usepackage{enumerate}
\usepackage[all,cmtip]{xy}
\usepackage{hyperref}

\usepackage{tikz}
\usetikzlibrary{decorations.pathreplacing,arrows}

\usepackage[margin=1in]{geometry}

\usepackage{setspace}

\allowdisplaybreaks

\usepackage{hyperref}

\newtheorem{theorem}{Theorem}[section]
\newtheorem{lemma}[theorem]{Lemma}

\theoremstyle {definition}
\newtheorem{remark}[theorem]{Remark}
\newtheorem{example}[theorem]{Example}

\numberwithin{theorem}{section}

\DeclareMathOperator{\area}{area} 
\newcommand{\R}{\mathbb{R}}
\DeclareMathOperator{\dist}{dist}

\title{A splitting theorem for scalar curvature}

\author{Otis Chodosh}
\address{Princeton University, Department of Mathematics, Fine Hall, Washington Road, Princeton, NJ 08544, United States}
\email{ochodosh@princeton.edu}

\author{Michael Eichmair}
\address{Faculty of Mathematics, University of Vienna, Oskar-Morgenstern-Platz 1, 1090 Vienna, Austria}
\email {michael.eichmair@univie.ac.at}

\author{Vlad Moraru}
\address{Department of Information Engineering, Computer Science and Mathematics, University of L'Aquila,
Via Vetoio - Coppito, 67100 l'Aquila, Italy}
\email{vlad.moraru@univaq.it}

\begin{document}

\begin{abstract}
We show that a Riemannian $3$-manifold with non-negative scalar curvature is flat if it contains an area-minimizing cylinder. This scalar-curvature analogue of the classical splitting theorem of J.~Cheeger and D.~Gromoll \cite{Cheeger-Gromoll:1971} has been conjectured by D.~Fischer-Colbrie and R.~Schoen \cite{Fischer-Colbrie-Schoen} and by M.~Cai and G.~Galloway \cite{Cai-Galloway:2000}. 
\end{abstract}

\maketitle

%%%%%%%%%%%%%%%%%%%%%%%%%%%%%%%%%%%%%%%%%
%%%%%%%%%%%%%%%%%%%%%%%%%%%%%%%%%%%%%%%%%
  
\section{Introduction}

Let $(M, g)$ be a connected, orientable, complete Riemannian $3$-manifold with non-negative scalar curvature. D.~Fischer-Colbrie and R.~Schoen show in \cite{Fischer-Colbrie-Schoen} that a connected, orientable, complete stable minimal immersion into $(M, g)$ is conformal to a plane, a sphere, a torus, or a cylinder. They conjecture that $(M, g)$ is flat if the immersion is conformal to the cylinder; cf.~Remark 4 in \cite{Fischer-Colbrie-Schoen}. M.~Cai and G.~Galloway point out a counterexample obtained from flattening standard $\R \times \mathbb{S}^2$ near $\R \times \{\text{great circle}\}$ in their concluding remark in \cite{Cai-Galloway:2000}. They ask if the conjecture holds under the additional assumption that the immersion be ``suitably" area-minimizing. In this paper, we prove the following result: 

\begin {theorem} \label{thm:main} 
Let $(M, g)$ be a connected, orientable, complete Riemannian $3$-manifold with non-negative scalar curvature. Assume that $(M, g)$ contains a properly embedded surface $S \subset M$ that is both homeomorphic to the cylinder and absolutely area-minimizing. Then $(M, g)$ is flat. In fact, a cover of $(M, g)$ is isometric to standard $\mathbb{S}^1 \times \R^2$ upon scaling. 
\end {theorem}

Note that this result is in satisfying analogy with the classical splitting theorem of J.~Cheeger and D.~Gromoll \cite{Cheeger-Gromoll:1971} in dimension $3$, where scalar curvature takes the place of Ricci curvature and where area-minimizing cylinders stand in for length-minimizing geodesic lines.\footnote{As the example of the doubled Schwarzschild manifold shows, the classical splitting theorem fails when we relax the assumption of non-negative Ricci curvature to non-negative scalar curvature.} 

Theorem \ref{thm:main} follows from the work of M.~Anderson and L.~Rodr\'iguez \cite{Anderson-Rodriguez:1989} when we impose the much stronger assumption of bounded, non-negative Ricci curvature. The strategy of M.~Anderson and L.~Rodr\'iguez \cite{Anderson-Rodriguez:1989} has been refined by G.~Liu \cite{Liu:2013} to classify complete, non-compact Riemannian $3$-manifolds with non-negative Ricci curvature. These ideas have been developed by the first- and second-named authors to establish the following \emph{scalar-curvature} rigidity result for asymptotically flat $3$-manifolds  which had been conjectured by R.~Schoen.

\begin {theorem} [\cite{mineffectivePMT}] \label{thm:Schoenrigidity} The only asymptotically flat Riemannian $3$-manifold with non-negative scalar curvature that admits a non-compact, area-minimizing boundary is flat $\R^3$. 
\end {theorem}

Our proof of Theorem \ref{thm:main} in this paper is a further development of these ideas. We now discuss challenges in the proof that are not present in the previously discussed works.

The goal is to construct a \textit{foliation} of $(M, g)$ by area-minimizing cylinders. The leaves of this foliation arise as limits of solutions of certain Plateau problems for least area. A major challenge we face here that has no substantial analogue in the proof of Theorem \ref{thm:Schoenrigidity} is the \textit{a priori} possible appearance of stable minimal planes or spheres (rather than cylinders or tori) in these limits. This scenario is addressed in Figure \ref{fig:no-planes} below. The ambient scalar curvature may well be positive along such surfaces, as the examples in Remark 3 of \cite{Fischer-Colbrie-Schoen} show. 

We need an approach that is sensitive to the topology of solutions of the Plateau problems in the construction. At the same time, we have to make sure that these solutions pass to limits in a reasonable way. Recall from e.g.~\cite{Rosenberg-Souam-Toubiana:2010} that stable minimal surfaces in Riemannian $3$-manifolds satisfy local curvature estimates that are independent of area bounds. In particular, sequences of such surfaces admit subsequential limits as pointed immersions. If each surface in the sequence is an area-minimizing boundary, then so is the limit. (If a small ambient ball intersects such a surface in two components, then these sheets have opposite orientation and are almost parallel. This scenario can be ruled out by a cut-and-paste argument.) However, limits of general area-minimizing surfaces can exhibit much greater complexity -- think of condensing closed geodesics on the torus and compare with Remark \ref{rem:stacks} below. 

As such, the use of solutions of Plateau problems in the class of \textit{all} (oriented) competitors risks the loss of local area bounds in the limit. On the other hand, the use of solutions in the class of boundaries risks the appearance of planes or spheres in the limit. These threats taken together force us to select the various classes of surfaces considered in the proof of Theorem \ref{thm:main} with great care.\\

We review  the notions of area-minimizing surfaces that are used in this paper in Appendix \ref{sec:am}. \\

Subsequent to the paper of M.~Cai and G.~Galloway \cite{Cai-Galloway:2000}, there have been many further works establishing scalar curvature rigidity results in the presence of \emph{compact} area-minimizing surfaces, including \cite{Cai:2002, Bray-Brendle-Eichmair-Neves:2010, Bray-Brendle-Neves:2010, Galloway:2011, MarquesNeves:min-max-rigidity-3mflds, MaximoNunes:hawking-rigidity, Nunes:2013,Ambrozio:free-bdry-rigidity, Micallef-Moraru:2015, Moraru:2016}. We anticipate that the techniques developed here lead to alternative proofs of these results. We plan to explore this possibility in forthcoming work. \\

Finally, we mention that parts of the strategy of M.~Anderson and L.~Rodr\'iguez \cite{Anderson-Rodriguez:1989} depend on the assumption of non-negative Ricci curvature in a subtle but essential way. In particular, we do not see how to carry over the crucial area estimate \cite[(1.5)]{Anderson-Rodriguez:1989} to the non-negative scalar curvature setting. \\

{\bf Acknowledgments: } We would like to thank Richard Bamler, Gregory Galloway, John Pardon, Leon Simon, and Brian White for  valuable discussions. The work of Otis Chodosh has been supported by EPSRC grant EP/K00865X/1, by the Oswald Veblen Fund, and by the NSF grant No.~1638352.  The work of Michael Eichmair has been supported by the START-Project Y963-N35 of the Austrian Science Fund (FWF). 

We dedicate this paper to Gregory Galloway on the occasion of his $70^{\text{th}}$ birthday.

%%%%%%%%%%%%%%%%%%%%%%%%%%%%%%%%%%%%%%%%
%%%%%%%%%%%%%%%%%%%%%%%%%%%%%%%%%%%%%%%%
%%%%%%%%%%%%%%%%%%%%%%%%%%%%%%%%%%%%%%%%

\section{Tools}

The following result is due to R.~Schoen and S.-T.~Yau \cite[Section 5]{Schoen-Yau:1979-Ann}, except for the assertion that $\varphi^* g$ is flat in the latter alternative, which is due to D.~Fischer-Colbrie and R.~Schoen \cite[Theorem 3]{Fischer-Colbrie-Schoen}.

\begin {lemma} \label{lem:fcsc} Let $(M, g)$ be a Riemannian $3$-manifold. Let $\varphi : S \to M$ be an orientable, complete, two-sided stable minimal immersion of a closed surface $S$ such that $R \circ \varphi \geq 0$, where $R$ is the scalar curvature of $(M, g)$. Then $S$ is topologically either a sphere or a torus. In the latter alternative, the immersion is totally geodesic, the induced metric $\varphi^* g$ on $S$ is flat, and $R \circ \varphi = 0$.
\end {lemma}

The following result due to D.~Fischer-Colbrie and R.~Schoen is part of Theorem 3 in \cite{Fischer-Colbrie-Schoen}. 

\begin {lemma}  \label{lem:fcs} Let $(M, g)$ be a Riemannian $3$-manifold. Let $\varphi : S \to M$ be an orientable, complete, non-compact, two-sided stable minimal immersion such that $R \circ \varphi \geq 0$ where $R$ is the scalar curvature of $(M, g)$. Then $S$ with the induced metric $\varphi^* g$ is conformal to either the plane or the cylinder. 
\end {lemma}

We refer to \cite{Schoen-Yau:1982, Miyaoka:1993, Berard-Castillon, Reiris} as well as Appendix C of \cite{mineffectivePMT} for discussions and proofs of the following refinement of the latter alternative in Lemma \ref{lem:fcs}. 

\begin {lemma} \label{lem:rigidcylinders} Let $(M, g)$ be a Riemannian $3$-manifold. Let $\varphi : S \to M$ be an orientable, complete, non-compact, two-sided stable minimal immersion such that $R \circ \varphi \geq 0$, where $R$ is the scalar curvature of $(M, g)$. If $S$ is a cylinder, then the immersion is totally geodesic, the induced metric $\varphi^* g$ is flat, and $R \circ \varphi = 0$.  
\end {lemma}

The following rigidity result is due to M.~Cai and G.~Galloway \cite[Theorem 2]{Cai-Galloway:2000}. 

\begin {lemma} \label{lem:Cai-Galloway} Let $(M, g)$ be a connected, complete Riemannian $3$-manifold with possibly empty weakly mean-convex boundary. We also assume that $(M, g)$ has non-negative scalar curvature. If $(M, g)$ contains a two-sided torus that has least area in its isotopy class, then $(M, g)$ is flat.
\end {lemma}

Finally, we will need the following lifting property for absolutely area-minimizing surfaces. 

\begin {lemma} \label{lem:topreduction} 
Let $(M, g)$ be an orientable, complete Riemannian $3$-manifold and let $S \subset M$ be a properly embedded, orientable, absolutely area-minimizing surface without boundary. There is a covering $p : \tilde M \to M$ with $p_* ( \pi_1(\tilde M)) = i_* ( \pi_1( S))$ where $i : S \to M$ is the inclusion map. The inclusion map lifts to a proper embedding $\tilde i : S \to \tilde M$ and $\tilde i_* : \pi_1(S) \to \pi_1(\tilde M)$ is surjective. Moreover, $\tilde S = \tilde i (S)  \subset \tilde M$ is absolutely area-minimizing in $(\tilde M, \tilde g)$ where $\tilde g = p^* g$.

\begin {proof} 
The asserted existence of covering and lift are standard, see Propositions 1.36 and 1.33 in \cite{Hatcher}. To see that $\tilde S$ is area-minimizing, note that $p$ is injective along $\tilde S$ and that every (properly embedded) competing surface for $\tilde S$ projects to a properly immersed, competing surface for $S$. Using cut-and-paste arguments with small changes in area, we obtain a properly embedded competitor downstairs.  
\end {proof}
\end {lemma}

%%%%%%%%%%%%%%%%%%%%%%%%%%%%%%%%%%%%%%%%
%%%%%%%%%%%%%%%%%%%%%%%%%%%%%%%%%%%%%%%%
%%%%%%%%%%%%%%%%%%%%%%%%%%%%%%%%%%%%%%%%

\section{Proof of Theorem \ref{thm:main}}

Let $S \subset M$ be a properly embedded cylinder that is absolutely area-minimizing in $(M,g)$. 

In view of Lemma \ref{lem:topreduction}, we may assume that the inclusion map $i : S \to M$ induces a surjection $i_* : \pi_1(S) \to \pi_1(M)$. It follows from Lemma \ref{lem:rigidcylinders} that $S \subset M$ is intrinsically flat. By scaling $(M, g)$ if necessary, we may assume that $S \subset M$ is isometric to standard $\mathbb{S}^1 \times \R$.

If $S \subset M$ is separating, then $M \setminus S$ has two components. We cut $M$ along $S$ and make a choice of component to obtain a connected, complete Riemannian $3$-manifold whose boundary is connected. If $S \subset M$ doesn't separate, then $M \setminus S$ is connected. We cut $M$ along $S$ to obtain a connected, complete Riemannian $3$-manifold whose boundary has exactly two components, of which we choose one. Either way, we denote the new Riemannian $3$-manifold by $(\hat M, \hat g)$ and the chosen component of its boundary by $\Sigma$. Note that $\Sigma \subset \hat M$ is isometric to $S$ and absolutely area-minimizing in $(\hat M, \hat g)$. We denote the closed curve on $\Sigma$ that corresponds to $\mathbb{S}^{1}\times \{0\}$ by $\gamma$. We write $\Sigma_{h}$ for the portion of $\Sigma$ corresponding to $\mathbb{S}^{1}\times [-h,h]$.

\begin{figure}[h]
\begin{tikzpicture}
	\begin{scope}[shift={(0,.17)}]
	\filldraw [opacity = .2] (0,1) circle (1.8);
	\filldraw [white] (0,1) circle (.6);
	\draw (0,1) circle (.6);
	\draw (0,1) circle (1.8);
	\filldraw (0,1) circle (1pt) node [shift={(-.0,.23)}] {$\scriptstyle c(2r)$};
	\draw [->] (0,1) -- node [shift = {(.1,.15)}] {$\scriptstyle r$} +(-40:.6);
	\draw [->] (0,1) -- node [shift = {(-.1,.2)}] {$\scriptstyle 3r$} +(185:1.8);
	\draw [->] (0,1) -- node [shift = {(.2,-.25)}] {$\scriptstyle 2r$} +(-92:1.15);
	\end{scope}
	\draw (0,0) to [bend left = 90, looseness = .5] (0,-2);
	\draw [dashed] (0,0) to [bend right = 90, looseness = .5] (0,-2);
	\begin{scope}[shift={(-6,0)}]
	\draw (0,0) to [bend left = 90, looseness = .5] (0,-2);
	\draw [dashed] (0,0) to [bend right = 90, looseness = .5] (0,-2);
	\end{scope}
	\begin{scope}[shift={(6,0)}]
	\draw (0,0) to [bend left = 90, looseness = .5] (0,-2);
	\draw [dashed] (0,0) to [bend right = 90, looseness = .5] (0,-2);
	\end{scope}
	\draw (-4,.48) to [bend left = 90, looseness = .5] (-4,-2.48);
	\draw [dashed] (-4,.48) to [bend right = 90, looseness = .5] (-4,-2.48);
	\draw (4,.48) to [bend left = 90, looseness = .5] (4,-2.48); 	\draw [dashed] (4,.48) to [bend right = 90, looseness = .5] (4,-2.48); 
	\draw (-6,0) to [bend left, looseness = .4] (6,0);
	\draw (-6,-2) to [bend right, looseness = .4] (6,-2);
	\filldraw [white, opacity =.3] (-4.5,-2) rectangle (-4.2,0);
	\filldraw [white, opacity =.3] (3.4,-2) rectangle (3.8,0);
	\draw (-7,0) -- (7,0);
	\draw (-7,-2) -- (7,-2);
	\node at (2,-1.7) {$\Sigma_{h} \cong \mathbb{S}^{1}\times [-h,h]$};
	\node at (-4,.8) {$\Sigma_{h}(r,t)$};	
	\node at (-4,2) {$\hat g(r,t) = \hat g$};
	\node at (4,2) {$R_{\hat g}\geq0$};
	\node at (0,2.2) {$R_{\hat g(r,t)} > 0$};
	\node at (-.5,-1.1) {$\gamma$};
\end{tikzpicture}
\caption{The scalar curvature of $\hat g(r,t)$ is positive in the shaded region.}
\label{fig:bumped-metric}
\end{figure}
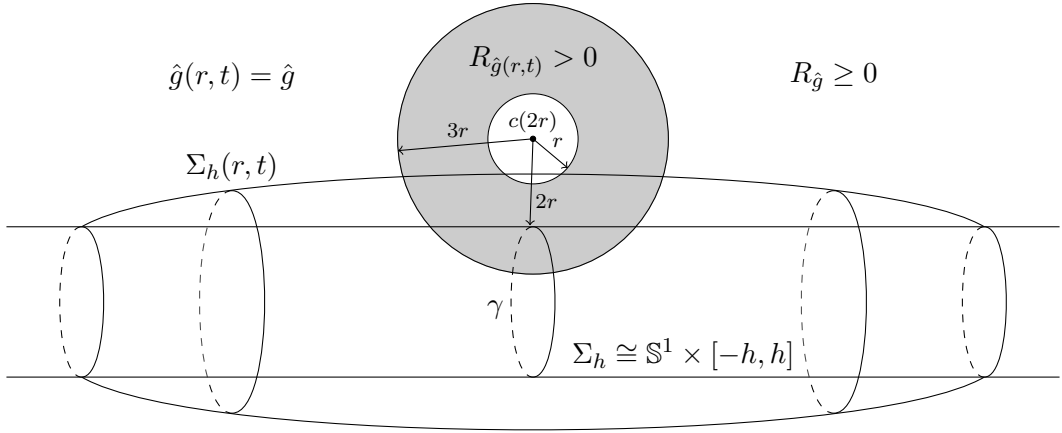

Fix a unit speed geodesic $c:[0,\varepsilon) \to \hat M$ with $c(0) \in \gamma$ and $\dot c (0) \perp \text{T}_{c(0)} \Sigma$. As in the proof of Theorem \ref{thm:Schoenrigidity} in Appendix J of \cite{mineffectivePMT}, we find a family of smooth Riemannian metrics $\{ \hat g (r, t)\}_{r, t \in (0, \varepsilon)}$ on $\hat M$ with the following properties (illustrated in Figure \ref{fig:bumped-metric}):
\begin {enumerate} [(i)] 
\item $  \hat g (r, t) \to \hat g$ in $C^3$ as $t,r \searrow 0$;
\item $ \hat  g (r, t) \to \hat g$ smoothly as $t \searrow 0$ for $r \in (0, \varepsilon)$ fixed;
\item $ \hat g (r, t) = \hat g$ on $\{x \in \hat M :  \dist_{\hat g} (x, c(2r)) \geq 3 r \}$;
\item $ \hat g (r, t) < \hat g$ as quadratic forms in $\{x \in \hat M : \dist_{\hat g} (x, c(2r)) < 3 r\}$;
\item $ \hat g (r, t)$ has positive scalar curvature in $\{x \in \hat M : r < \dist_{\hat g}(x, c(2r)) < 3 r\}$;
\item $\hat M$ is weakly mean-convex with respect to $\hat g (r, t)$. 
\end {enumerate}
Fix $h > 1$. Let $B_{h}$ denote a pre-compact open set with smooth boundary in $\hat M$ and such that $\{x \in \hat M : \dist_{\hat g}(x,\Sigma_{h}) < 2h\} \subset B_h$. We modify the metric $\hat g (r,t)$ near the boundary of $B_{h}$ to $\hat g(r,t,h)$ so $B_{h}$ is weakly mean-convex with respect to $\hat g(r,t,h)$ and 
\[
(1-\delta) \hat g(r,t) \leq \hat g (r,t,h) \leq (1+\delta) \hat g (r,t)
\]
where $\delta \in (0, 1)$ is chosen to satisfy \eqref{choicedelta} below. 
Among all compact, oriented surfaces in $B_h$ with boundary $\partial \Sigma_{h}$ that bound an open subset of $\hat M$, there is one whose area with respect to $\hat g (r, t, h)$ is least. Choose one such area-minimizing surface and denote it by $\Sigma_{h}(r,t)$. 
We claim that $\Sigma_{h}(r,t)$ intersects $\{x \in \hat M : \dist_{\hat g}(x,c(2r)) < 3r\}$. For if not, we have that $\hat g= \hat g(r,t)$ along $\Sigma_{h}(r,t)$ and we may compute, as in \cite[(12)]{mineffectivePMT}, that
\begin{align*}
0 < \area_{\hat g}(\Sigma_{h}) - \area_{\hat g(r,t)}(\Sigma_{h}) & \leq \area_{\hat g}(\Sigma_{h}(r,t)) - \area_{\hat g(r,t)}(\Sigma_{h})\\
& = \area_{\hat g(r,t)}(\Sigma_{h}(r,t)) - \area_{\hat g(r,t)}(\Sigma_{h}) \\
& \leq \frac{1}{1 - \delta} \area_{\hat g(r,t,h)}(\Sigma_{h}(r,t)) - \area_{\hat g(r,t)}(\Sigma_{h}) \\
& \leq  \frac{1}{1 - \delta} \area_{\hat g(r,t,h)}(\Sigma_{h}) - \area_{\hat g(r,t)}(\Sigma_{h}) \\
& \leq \frac{1+\delta}{1-\delta} \area_{\hat g(r,t)}(\Sigma_{h}) - \area_{\hat g(r,t)}(\Sigma_{h}) \\
& = \frac{2\delta}{1-\delta} \area_{\hat g(r,t)}(\Sigma_{h}).
\end{align*}
This is a contradiction if we choose $\delta =  \delta (r, t, h) >0$ with
\begin{align} \label{choicedelta}
\frac{2\delta}{1-\delta} < \frac{\area_{\hat g}(\Sigma_{h}) - \area_{\hat g(r,t)}(\Sigma_{h})}{\area_{\hat g(r,t)}(\Sigma_{h})}.
\end{align}

A comparison with small geodesic spheres gives local area bounds for the surfaces $\Sigma_h(r, t)$ that are independent of all parameters; cf.~Remark \ref{rem:stacks}. Using standard results in geometric measure theory, we may pass these surfaces to a subsequential limit as $h\to\infty$ to obtain a properly embedded surface $\Sigma (r, t)$. Note that $\Sigma (r, t)$ is a boundary in $\hat M$. It follows from the construction that this boundary is homologically* area-minimizing\footnote{This non-standard area-minimizing property is discussed in Appendix \ref{sec:am}. Its full strength is needed in the cut-and-paste argument used to rule out planes and spheres below.} with respect to $g(r, t)$. 

When $r, t > 0$ are sufficiently small, the surface $\Sigma(r,t)$ contains a closed curve $\gamma (r,t)$ that intersects $\{x\in \hat M : \dist_{\hat g}(x,c(2r)) \leq 3r\}$ and which is close to $\gamma$. Indeed, there is a small ball with center at $c(0) \in \gamma$ where a piece of $\Sigma (r,t)$ appears as a graph above $\text{T}_{c(0)} \Sigma$. The geometric Harnack principle allows us to continue this piece of $\Sigma(r, t)$ into a ribbon as we traverse $\gamma$. The ribbon must close up as we travel around $\gamma$, as otherwise there would be two nearby sheets -- a scenario that contradicts the area-minimizing property of $\Sigma(r, t)$. The component $\hat \Sigma(r,t)$ of $\Sigma(r,t)$ that contains this curve $\gamma(r,t)$ converges to $\Sigma$ as $r, t \searrow 0$. Observe that $\hat \Sigma (r, t)$ is a boundary in $\hat M$ that is homologically* area-minimizing with respect to $\hat g(r, t)$; cf.~Lemma 33.4 in \cite{GMT}. 

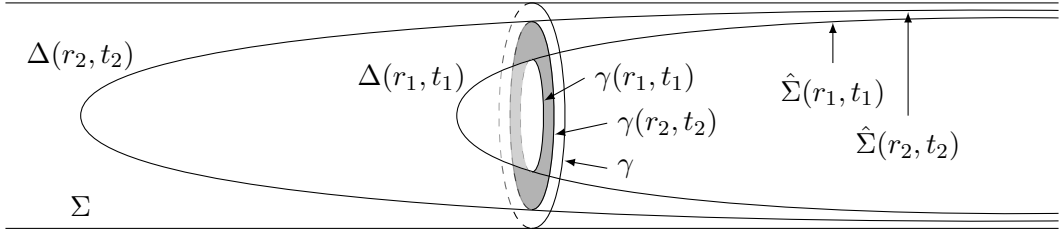
\begin{figure}[h!]
\begin{tikzpicture}
	\draw (-7,0) -- (7,0);
	\draw (-7,-3) -- (7,-3);
	\draw [dashed] (0,-.25) to [bend right = 90, looseness = .4] (0,-2.75);
	\draw [dashed] (0,-.75) to [bend right = 90, looseness = .35] (0,-2.25);	
	\draw [dashed] (0,0) to [bend right = 90, looseness = .5] (0,-3);
	\filldraw [white, opacity = .4] plot [smooth, tension = 2] coordinates {(7,-.1) (-6,-1.5) (7,-2.9)};
	\filldraw [opacity = .3] (0,-.25) to [bend left = 90, looseness = .4] (0,-2.75);
	\filldraw [white] (-.01,-.75) to [bend left = 90, looseness = .35] (-.01,-2.25);		
	\filldraw [opacity = .3] (0,-.25) to [bend right = 90, looseness = .4] (0,-2.75);
	\filldraw [white] (.01,-.75) to [bend right = 90, looseness = .35] (.01,-2.25);	
	\draw plot [smooth, tension =2]	coordinates {(7,-.1) (-6,-1.5) (7,-2.9)};
	\begin{scope}
	\clip (-1.1,-3) rectangle (0,0);
	\filldraw [white, opacity =.4] plot [smooth, tension = 2] coordinates {(7,-.2) (-1,-1.5) (7,-2.8)};
	\end{scope}
	\draw plot [smooth, tension = 2] coordinates {(7,-.2) (-1,-1.5) (7,-2.8)};
	\draw (0,0) to [bend left = 90, looseness = .5] (0,-3);
	\draw (0,-.75) to [bend left = 90, looseness = .35] (0,-2.25);
	\draw (0,-.25) to [bend left = 90, looseness = .4] (0,-2.75);
	\node at (-6,-.7) {$\Delta(r_{2},t_{2})$};
	\node at (-1.6,-1) {$\Delta(r_{1},t_{1})$};
	\draw [-latex] (.7,-1) node [right] {$\gamma(r_{1},t_{1})$} -- (.16,-1.3);
	\draw [-latex] (1,-1.6) node [right] {$\gamma(r_{2},t_{2})$} -- (.3,-1.75);
	\draw [-latex] (1,-2.2) node [right] {$\gamma$} -- (.43,-2.1);	
	\draw [-latex] (4,-.8) node [below] {$\hat\Sigma(r_{1},t_{1})$} -- (4,-.26);
	\draw [-latex] (5,-1.5) node [below] {$\hat\Sigma(r_{2},t_{2})$} -- (5,-.12);
	\node at (-6,-2.7) {$\Sigma$};
\end{tikzpicture}
\caption{The possibility that $\hat\Sigma(r,t)$ is a plane for $r,t>0$ small enough can be ruled out by comparing the areas of the disks $\Delta(r,t)$, or by incompressibility of $\Sigma$.}
\label{fig:no-planes}
\end{figure}

We claim that $\hat \Sigma (r, t)$ is \textit{neither} a plane \textit{nor} a sphere if $r, t > 0$ are small enough. Suppose otherwise. In this case, $\gamma(r,t)$ bounds an embedded disk in $\hat M$. Consequently, so does $\gamma$. In particular, inclusion induces the trivial map $\pi_{1}(\Sigma)\to \pi_{1}(\hat M)$. By Lemma \ref{lemm:key-top-lemma} below, every connected, closed surface in $\hat M$ is separating. If $\hat \Sigma (r, t)$ is a sphere, it bounds a possibly unbounded region in $\hat M$. This contradicts its homologically* area-minimizing property. Assume that $\hat \Sigma(r,t)$ is a plane (see Figure \ref{fig:no-planes}). Let $\Delta(r,t) \subset \hat \Sigma(r,t)$ be the disk bounded by $\gamma (r, t)$. Since $\hat \Sigma(r,t)$ converges to $\Sigma$ as $r,t \searrow 0$, we have that $\area_{g}(\Delta(r,t)) \to \infty$ as $r,t\searrow 0$. Choose $0 < r_{2}\ll r_{1}$ and $0< t_{2} \ll t_{1}$ such that $\hat \Sigma (r_1, t_1)$ and $\hat \Sigma (r_2, t_2)$ are planes. Note that $\Delta(r_{2},t_{2})$, $\Delta(r_{1},t_{1})$, and a small neck connecting $\gamma(r_{1},t_{1})$ and $ \gamma(r_{2},t_{2})$ bound an open set. This contradicts the homologically* area-minimizing property of $\hat \Sigma(r_{1},t_{1})$. 

We see from (v), Lemma \ref{lem:fcsc}, Lemma \ref{lem:fcs}, and Lemma \ref{lem:rigidcylinders} that, for $r, t > 0$ small enough, $\hat \Sigma(r,t)$ either intersects $\{x \in \hat M : \dist_{\hat g}(x,c(2r)) \leq r\}$, or $\hat \Sigma(r,t)$ intersects $\{x \in \hat M : \dist_{\hat g}(x,c(2r)) = 3r\}$ but not $\{x \in \hat M : \dist_{\hat g}(x,c(2r)) < 3r\}$. Fix $r > 0$ small and pass to a geometric subsequential limit as $t \searrow 0$. We obtain a properly embedded boundary $\Sigma(r) \subset \hat M$ which has a connected component $\hat \Sigma (r)$ that intersects $\{x\in \hat M : \dist_{\hat g}(x,c(2r)) \leq 3r\}$ and which is disjoint from $\Sigma$. As before, we see that $\hat \Sigma (r)$ is a properly embedded boundary in $\hat M$ that is homologically* area-minimizing with respect to $\hat g$. Clearly, $\hat \Sigma (r)$ is disjoint from $\Sigma$ and contains a closed embedded curve $\gamma(r)$ close to $\gamma$. The argument of the preceding paragraph shows that when $r>0$ is small enough, the surface $\hat \Sigma(r)$ is diffeomorphic to either a torus or a cylinder.

If $\hat \Sigma (r)$ is a torus, then $(\hat M, \hat g)$ is isometric to either standard $\mathbb{S}^1 \times \R \times [0, \infty)$ or standard $\mathbb{S}^1 \times \R \times [0, a]$ for some $a > 0$ by Lemma \ref{lem:Cai-Galloway}. We may thus assume that $\hat \Sigma(r)$ is cylindrical for all $r > 0$ small. By Lemma \ref{lem:rigidcylinders}, $\hat \Sigma (r) \subset \hat M$ is intrinsically flat, totally geodesic, and the ambient Ricci tensor evaluated in the normal direction vanishes along $\hat \Sigma (r)$. Note that $\hat \Sigma(r)$ converges to $\Sigma$ as $r \searrow 0$. We now argue exactly as in the proof of \cite[Theorem 1.6]{mineffectivePMT} to show that the ambient Riemann tensor vanishes along $\Sigma$. 

We may repeat the above argument starting with any of the surfaces $\hat \Sigma(r)$ for $r > 0$ sufficiently small, using that they are homologically* area-minimizing in $(\hat M, \hat g)$.\footnote{The advantage of the \emph{absolutely} area-minimizing property in the above proof is that it lifts to covers. We only used this at the very beginning, when we applied Lemma \ref{lem:topreduction}.} A continuity argument then gives that $(\hat M, \hat g)$ is either standard $\mathbb{S}^1 \times \R \times [0, \infty)$ or standard $\mathbb{S}^1 \times \R \times [0, a]$ for some $a > 0$. \qed 

\begin{lemma} \label{lemm:key-top-lemma} 
Assumptions and notation as in the proof of Theorem \ref{thm:main} above. If the inclusion $\Sigma \subset \hat M$ induces a trivial map $\pi_1(\Sigma) \to \pi_1(\hat M)$, then every connected, closed surface in $\hat M$ is separating. 
\begin{proof} 
It follows from the hypothesis and the construction of $\hat M$ from $M \setminus S$ that $M$ is simply connected. (Here we are use that the map $\pi_1(S) \to \pi_1(M)$ induced by the inclusion $S \subset M$ is surjective.) Standard intersection theory gives that every connected, closed embedded surface $N \subset M$ separates $M$. Indeed, if we assume that $N \subset M$ is not separating, there is a closed embedded curve in $M$ that intersects $N$ transversely and exactly once. Such a curve cannot be homotopically trivial. Assume now that $\hat N \subset \hat M$ is a connected and closed surface disjoint from the boundary of $\hat M$. It corresponds to a connected and closed surface $N \subset M$ that is disjoint from $S$. Using that $N$ separates $M$ and the construction of $\hat M$ from $M$, we conclude that $\hat N$ separates. 
\end{proof}
\end{lemma}

A variation of this proof of Theorem \ref{thm:main} gives the following result. 

\begin{theorem} \label{thm:incomp-variant} Let $(M, g)$ be a connected, orientable, complete Riemannian $3$-manifold with non-negative scalar curvature. Assume that $(M, g)$ contains a properly embedded, incompressible\footnote{More precisely, we require that every loop $\gamma \subset S$ that bounds an embedded disk $\Delta \subset M$ is contractible in $S$.} surface $S \subset M$ that is homeomorphic to the cylinder and homologically area-minimizing. Then $(M, g)$ is flat. In fact, a cover of $(M, g)$ is isometric to standard $\mathbb{S}^1 \times \R^2$ upon scaling. 
\begin{proof}
We follow the proof of Theorem \ref{thm:main} above, except for the following changes: 
\begin {enumerate} [(a)]
\item We do not pass to a covering. 
\item We work with surfaces $\Sigma_h(r, t)$ with $\partial \Sigma_h(r, t) =\partial \Sigma_h$ that have least area with respect to $\hat g(r, t, h)$ and which together with $\Sigma_h$ bound an open set in $B_h$. 
\item To argue that $\hat \Sigma (r, t)$ is neither a plane nor a sphere, we use that $\Sigma$ is incompressible.\footnote{It is not clear that the comparison surfaces we use in the cut-and-paste argument in the proof of Theorem \ref{thm:main} above are homologous to the original surface. A priori, they may differ by an \textit{unbounded} open set.}\qedhere
\end {enumerate}
\end{proof}
\end {theorem}

\begin{remark}
We do not know if Theorem \ref{thm:incomp-variant} holds if we drop the assumption that $S$ be incompressible. 
\end{remark}

%%%%%%%%%%%%%%%%%%%%%%%%%%%%%%%%%%%%%%%%
%%%%%%%%%%%%%%%%%%%%%%%%%%%%%%%%%%%%%%%%
%%%%%%%%%%%%%%%%%%%%%%%%%%%%%%%%%%%%%%%%

\appendix

\section {Notions of area-minimizing surfaces} \label{sec:am}

Let $(M, g)$ be an orientable Riemannian manifold -- possibly with boundary. 

Let $\Sigma \subset M$ be an oriented, properly embedded hypersurface. 

Recall that $\Sigma \subset M$ is \emph{absolutely area-minimizing} in $(M, g)$ if for every $U \subset M$ open with compact closure, we have that
\[
\text{area} (U \cap \Sigma) \leq \text{area} (U \cap \tilde \Sigma)
\]
whenever $\tilde \Sigma \subset M$ is an oriented, properly embedded hypersurface with $\partial \tilde \Sigma = \partial \Sigma$ (matching orientations) and
\[
\tilde \Sigma \setminus U = \Sigma \setminus U.
\]

Recall that $\Sigma \subset M$ is \emph{homologically area-minimizing} in $(M, g)$ if for every  $U \subset M$ open with compact closure, we have that
\[
\text{area} (U \cap \Sigma) \leq \text{area} (U \cap \tilde \Sigma)
\]
whenever $\tilde \Sigma \subset M$ is an oriented, properly embedded hypersurface such that 
\[
\tilde \Sigma =  \Sigma + \partial \Omega_1 + \ldots + \partial \Omega_N 
\]
in the sense of Stokes' theorem, where $\Omega_1, \ldots ,\Omega_N \subset M$ are compact top-dimensional submanifolds with $\Omega_i \subset U$.  

\begin {remark} \label {rem:stacks} 
There is a standard a priori area bound for homologically area-minimizing boundaries; cf.~\cite[\S 37.2]{GMT}. The properties of such boundaries are preserved under convergence. There is no such a priori bound for general area-minimizing surfaces. Indeed, every stack of finitely many parallel planes $\R^n \times \{z_i\}$ with standard orientation and where $z_1, \ldots, z_m \in \R$ is absolutely area-minimizing in $\R^{n+1}$. 
\end {remark}

Finally, we discuss a non-standard notion that plays a pivotal role in the proof of Theorem \ref{thm:main}.

We say that $\Sigma \subset M$ is \emph{homologically* area-minimizing} in $(M, g)$ if for every  $U \subset M$ open with compact closure, we have that
\[
\text{area} (U \cap \Sigma) \leq \text{area} (U \cap \tilde \Sigma)
\]
whenever $\tilde \Sigma \subset M$ is an oriented, properly embedded hypersurface such that 
\[
\tilde \Sigma =  \Sigma + \partial \Omega_1 + \ldots + \partial \Omega_N 
\]
in the sense of Stokes' theorem, where $\Omega_1, \ldots ,\Omega_N \subset M$ are properly embedded top-dimensional submanifolds with $\partial \Omega_i \subset U$. \textit{The point is that we do \emph{not} require that $\Omega_i \subset U$ or even that $\Omega_i$ is bounded here.} 

\begin {example} Consider the embedded curves 
\begin {align*}
\gamma_1 &= \{ (e^{i \theta}, 0) : 0 \leq \theta \leq \pi/2\}  \\
\gamma_2 &= \{ (e^{i \theta}, 0 ) : \pi/2 \leq \theta \leq 2 \pi\}
\end {align*}
in the standard cylinder $\mathbb{S}^1 \times \R$. $\gamma_1$ is absolutely length-minimizing, $\gamma_2$ is not. Both $\gamma_1$ and $\gamma_2$ are homologically length-minimizing. $\gamma_1$ is homologically* length-minimizing, $\gamma_2$ is not.
\end {example}

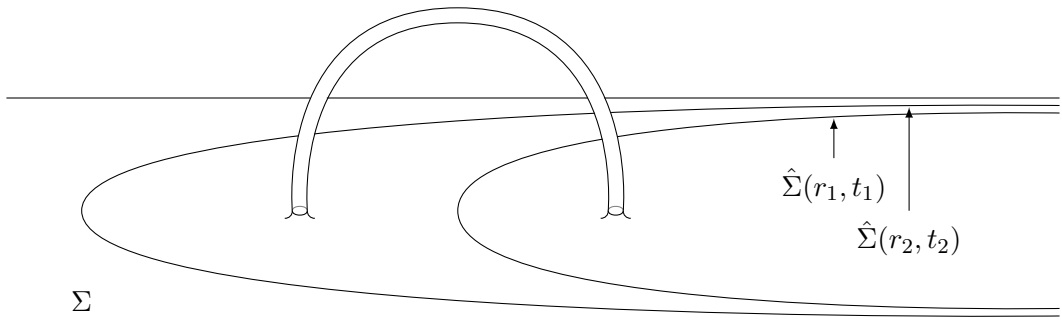
\begin{figure}[h!]
\vspace{2pt}
\begin{tikzpicture}
	\draw (-7,0) -- (7,0);
	\draw (-7,-3) -- (7,-3);
	\filldraw [white, opacity = .3] plot [smooth, tension = 2] coordinates {(7,-.1) (-6,-1.5) (7,-2.9)};	\draw plot [smooth, tension = 2] coordinates {(7,-.1) (-6,-1.5) (7,-2.9)};	
	\filldraw [white, opacity =.3] plot [smooth, tension = 2] coordinates {(7,-.2) (-1,-1.5) (7,-2.8)};
	\draw plot [smooth, tension = 2] coordinates {(7,-.2) (-1,-1.5) (7,-2.8)};
	\draw [-latex] (4,-.8) node [below] {$\hat\Sigma(r_{1},t_{1})$} -- (4,-.27);
	\draw [-latex] (5,-1.5) node [below] {$\hat\Sigma(r_{2},t_{2})$} -- (5,-.12);
	\node at (-6,-2.7) {$\Sigma$};
	
	\filldraw [white] (-3.01,-.49) circle (2.7 pt);
	\filldraw [white] (-2.85,0) circle (2.9 pt);
	\filldraw [white] (.84,0) circle (2.9 pt);	
	\filldraw [white] (.93,-.2) circle (2.9 pt);	
	\filldraw [white] (1.02,-.51) circle (2.7 pt);	
	
	\draw plot [smooth, tension = 2] coordinates {(-3,-1.5) (-1,1) (1,-1.5)};
	\draw plot [smooth, tension = 2] coordinates {(-3.2,-1.5) (-1,1.2) (1.2,-1.5)};
	
	\draw (-3,-1.5) to [bend left = 90] (-3.2,-1.5);
	\draw [opacity=.5] (-3,-1.5) to [bend right = 90] (-3.2,-1.5);
	
	\draw  [opacity=.5]  (1,-1.5) to [bend left = 90] (1.2,-1.5);
	\draw  (1,-1.5) to [bend right = 90] (1.2,-1.5);
	
	\draw (-3.2,-1.5) to [bend left = 50] (-3.3,-1.6);
	\draw (-3,-1.5) to [bend right = 50] (-2.9,-1.6);
	
	\draw (1,-1.5) to [bend left = 50] (.9,-1.6);
	\draw (1.2,-1.5) to [bend right = 50] (1.3,-1.6);
\end{tikzpicture}
\caption{This figure shows a hypothetical scenario if we had not passed to the cover (using Lemma \ref{lem:topreduction}). We want to compare the areas of $\hat\Sigma(r_{1},t_{1})$ and $\hat \Sigma(r_{2},t_{2})$ in a bounded set, using their respective minimizing properties. However, due to the presence of a \textit{neck}, the surfaces are neither homologically nor homologically* related.}
\label{fig:necks}
\end{figure}

\begin{figure}[h!]
\begin{tikzpicture}
	\draw (-7,0) -- (7,0);
	\draw (-7,-3) -- (7,-3);
	\begin{scope}
		\clip  (-6.1,-3) rectangle (6,0);
		\filldraw [black, opacity = .2] plot [smooth, tension = 2] coordinates {(7,-.1) (-6,-1.5) (7,-2.9)};
	\end{scope}
	\draw plot [smooth, tension = 2] coordinates {(7,-.1) (-6,-1.5) (7,-2.9)};	
	\filldraw [white] plot [smooth, tension = 2] coordinates {(7,-.2) (-1,-1.5) (7,-2.8)};
	\draw plot [smooth, tension = 2] coordinates {(7,-.2) (-1,-1.5) (7,-2.8)};
	\draw [-latex] (4,-.8) node [below] {$\hat\Sigma(r_{1},t_{1})$} -- (4,-.27);
	\draw [-latex] (5,-1.5) node [below] {$\hat\Sigma(r_{2},t_{2})$} -- (5,-.1);
	\node at (-6,-2.7) {$\Sigma$};
	
	\filldraw [white] (-3.2,-1.5) rectangle (-3,1);
%	\filldraw [black, opacity=.2] (-3.2,-1.5) rectangle (-3,1);
	\filldraw [black, opacity=.2] plot coordinates { (-3.2,-1.5) (-3.2,1.01) (-3.28,1.1) (-2.92,1.1) (-3,1.01) (-3,-1.5)};
	\filldraw [white] (-3.2,1) to [bend right = 50] (-3.3,1.1);
	\filldraw [white] (-3,1) to [bend left = 50] (-2.9,1.1);

	\draw (-3,-1.5) -- (-3,1);
	\draw (-3.2,-1.5) -- (-3.2,1);
	
	\draw (-3,-1.5) to [bend left = 90] (-3.2,-1.5);
	\draw [opacity=.5] (-3,-1.5) to [bend right = 90] (-3.2,-1.5);

	\filldraw [white] (1,-1.5) rectangle (1.2,1);
	\draw (1,-1.5) -- (1,1);
	\draw (1.2,-1.5) -- (1.2,1);
	\draw  [opacity=.5]  (1,-1.5) to [bend left = 90] (1.2,-1.5);
	\draw  (1,-1.5) to [bend right = 90] (1.2,-1.5);

	\draw  [opacity=.5]  (1,.8) to [bend left = 90] (1.2,.8);
	\draw  (1,.8) to [bend right = 90] (1.2,.8);

	\draw  [opacity=.5]  (-3.2,.8) to [bend left = 90] (-3,.8);
	\draw  (-3.2,.8) to [bend right = 90] (-3,.8);

	\node at (-3.1,1.4) {$\vdots$};	
	\node at (1.1,1.4) {$\vdots$};
	
	\draw (-3.2,-1.5) to [bend left = 50] (-3.3,-1.6);
	\draw (-3,-1.5) to [bend right = 50] (-2.9,-1.6);
	
	\draw (1,-1.5) to [bend left = 50] (.9,-1.6);
	\draw (1.2,-1.5) to [bend right = 50] (1.3,-1.6);
	
	\draw (-3.2,1) to [bend right = 50] (-3.3,1.1);
	\draw (-3,1) to [bend left = 50] (-2.9,1.1);
	
	\draw (1,1) to [bend right = 50] (.9,1.1);
	\draw (1.2,1) to [bend left = 50] (1.3,1.1);	
\end{tikzpicture}
\caption{A depiction of the situation in Figure \ref{fig:necks} after passing to the appropriate cover (using Lemma \ref{lem:topreduction}). Now, compact pieces of $\hat \Sigma(r_{1},t_{1})$ and $\hat \Sigma(r_{2},t_{2})$ are homologically* but not homologically comparable. This is because they differ by an open set (shaded grey) which is necessarily unbounded. }
\label{fig:covering.necks}
\end{figure}
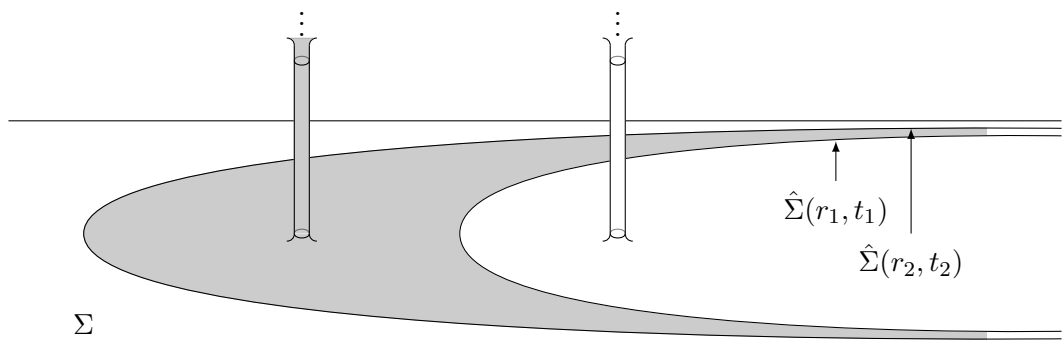

We conclude by describing how this non-standard area-minimizing property is tailored to a delicate aspect in the proof of Theorem \ref{thm:main}.  As we have just discussed, we cannot expect local area bounds for sequences of absolutely area-minimizing surfaces. For this reason, we want to fix the homology class of the surfaces we work with, or at least keep it in check. If we work with homological area-minimizers, it is not clear how to rule out minimizing planes or spheres in the proof of Theorem \ref{thm:main} by cut-and-paste arguments. Figure \ref{fig:necks} shows a scenario where the planes we worry about are \textit{not} homologically related. To deal with this scenario, we pass to the cover using Lemma \ref{lem:topreduction} at the beginning of the proof. The situation after passing to the cover is shown in Figure \ref{fig:covering.necks}. Now, the surfaces whose areas we want to compare bound a \textit{non-compact} region. They are homologically* related, but not homologically related.

\bibliography{bib} 

\providecommand{\bysame}{\leavevmode\hbox to3em{\hrulefill}\thinspace}
\providecommand{\MR}{\relax\ifhmode\unskip\space\fi MR }
\providecommand{\MRhref}[2]{%
  \href{http://www.ams.org/mathscinet-getitem?mr=#1}{#2}
}
\providecommand{\href}[2]{#2}
\begin{thebibliography}{10}

\bibitem{Ambrozio:free-bdry-rigidity}
Lucas Ambrozio, \emph{Rigidity of area-minimizing free boundary surfaces in
  mean convex three-manifolds}, J. Geom. Anal. \textbf{25} (2015), no.~2,
  1001--1017. \MR{3319958}

\bibitem{Anderson-Rodriguez:1989}
Michael Anderson and Lucio Rodr{\'{\i}}guez, \emph{Minimal surfaces and
  {$3$}-manifolds of nonnegative {R}icci curvature}, Math. Ann. \textbf{284}
  (1989), no.~3, 461--475. \MR{1001714}

\bibitem{Berard-Castillon}
Pierre B{\'e}rard and Philippe Castillon, \emph{Inverse spectral positivity for
  surfaces}, Rev. Mat. Iberoam. \textbf{30} (2014), no.~4, 1237--1264.
  \MR{3293432}

\bibitem{Bray-Brendle-Eichmair-Neves:2010}
H.~Bray, S.~Brendle, M.~Eichmair, and A.~Neves, \emph{Area-minimizing
  projective planes in 3-manifolds}, Comm. Pure Appl. Math. \textbf{63} (2010),
  no.~9, 1237--1247. \MR{2675487}

\bibitem{Bray-Brendle-Neves:2010}
Hubert Bray, Simon Brendle, and Andr{\'e} Neves, \emph{Rigidity of
  area-minimizing two-spheres in three-manifolds}, Comm. Anal. Geom.
  \textbf{18} (2010), no.~4, 821--830. \MR{2765731}

\bibitem{Cai:2002}
Mingliang Cai, \emph{Volume minimizing hypersurfaces in manifolds of
  nonnegative scalar curvature}, Minimal surfaces, geometric analysis and
  symplectic geometry ({B}altimore, {MD}, 1999), Adv. Stud. Pure Math.,
  vol.~34, Math. Soc. Japan, Tokyo, 2002, pp.~1--7. \MR{1925731}

\bibitem{Cai-Galloway:2000}
Mingliang Cai and Gregory Galloway, \emph{Rigidity of area minimizing tori in
  {$3$}-manifolds of nonnegative scalar curvature}, Comm. Anal. Geom.
  \textbf{8} (2000), no.~3, 565--573. \MR{1775139}

\bibitem{mineffectivePMT}
Alessandro Carlotto, Otis Chodosh, and Michael Eichmair, \emph{Effective
  versions of the positive mass theorem}, Invent. Math. \textbf{206} (2016),
  no.~3, 975--1016. \MR{3573977}

\bibitem{Cheeger-Gromoll:1971}
Jeff Cheeger and Detlef Gromoll, \emph{The splitting theorem for manifolds of
  nonnegative {R}icci curvature}, J. Differential Geometry \textbf{6}
  (1971/72), 119--128. \MR{0303460}

\bibitem{Fischer-Colbrie-Schoen}
Doris Fischer-Colbrie and Richard Schoen, \emph{The structure of complete
  stable minimal surfaces in {$3$}-manifolds of nonnegative scalar curvature},
  Comm. Pure Appl. Math. \textbf{33} (1980), no.~2, 199--211. \MR{562550}

\bibitem{Galloway:2011}
Gregory~J. Galloway, \emph{Stability and rigidity of extremal surfaces in
  {R}iemannian geometry and general relativity}, Surveys in geometric analysis
  and relativity, Adv. Lect. Math. (ALM), vol.~20, Int. Press, Somerville, MA,
  2011, pp.~221--239. \MR{2906927}

\bibitem{Hatcher}
Allen Hatcher, \emph{Algebraic topology}, Cambridge University Press,
  Cambridge, 2002. \MR{1867354}

\bibitem{Liu:2013}
Gang Liu, \emph{{$3$}-manifolds with nonnegative {R}icci curvature}, Invent.
  Math. \textbf{193} (2013), no.~2, 367--375. \MR{3090181}

\bibitem{MarquesNeves:min-max-rigidity-3mflds}
Fernando Marques and Andr{\'e} Neves, \emph{Rigidity of min-max minimal spheres
  in three-manifolds}, Duke Math. J. \textbf{161} (2012), no.~14, 2725--2752.
  \MR{2993139}

\bibitem{MaximoNunes:hawking-rigidity}
Davi M{\'a}ximo and Ivaldo Nunes, \emph{Hawking mass and local rigidity of
  minimal two-spheres in three-manifolds}, Comm. Anal. Geom. \textbf{21}
  (2013), no.~2, 409--432. \MR{3043752}

\bibitem{Micallef-Moraru:2015}
Mario Micallef and Vlad Moraru, \emph{Splitting of {$3$}-manifolds and rigidity
  of area-minimising surfaces}, Proc. Amer. Math. Soc. \textbf{143} (2015),
  no.~7, 2865--2872. \MR{3336611}

\bibitem{Miyaoka:1993}
Reiko Miyaoka, \emph{{$L^2$} harmonic {$1$}-forms on a complete stable minimal
  hypersurface}, Geometry and global analysis ({S}endai, 1993), Tohoku Univ.,
  Sendai, 1993, pp.~289--293. \MR{1361194}

\bibitem{Moraru:2016}
Vlad Moraru, \emph{On area comparison and rigidity involving the scalar
  curvature}, J. Geom. Anal. \textbf{26} (2016), no.~1, 294--312. \MR{3441515}

\bibitem{Nunes:2013}
Ivaldo Nunes, \emph{Rigidity of area-minimizing hyperbolic surfaces in
  three-manifolds}, J. Geom. Anal. \textbf{23} (2013), no.~3, 1290--1302.

\bibitem{Reiris}
Martin Reiris, \emph{Geometric relations of stable minimal surfaces and
  applications}, preprint, \url{http://arxiv.org/abs/1002.3274} (2010).

\bibitem{Rosenberg-Souam-Toubiana:2010}
Harold Rosenberg, Rabah Souam, and Eric Toubiana, \emph{General curvature
  estimates for stable {H}-surfaces in {$3$}-manifolds and applications}, J.
  Differential Geom. \textbf{84} (2010), no.~3, 623--648. \MR{2669367}

\bibitem{Schoen-Yau:1979-Ann}
Richard Schoen and Shing-Tung Yau, \emph{Existence of incompressible minimal
  surfaces and the topology of three-dimensional manifolds with nonnegative
  scalar curvature}, Ann. of Math. (2) \textbf{110} (1979), no.~1, 127--142.
  \MR{0541332}

\bibitem{Schoen-Yau:1982}
\bysame, \emph{Complete three-dimensional manifolds with positive {R}icci
  curvature and scalar curvature}, Seminar on {D}ifferential {G}eometry, Ann.
  of Math. Stud., vol. 102, Princeton Univ. Press, Princeton, N.J., 1982,
  pp.~209--228. \MR{645740}

\bibitem{GMT}
Leon Simon, \emph{Lectures on geometric measure theory}, Proceedings of the
  Centre for Mathematical Analysis, Australian National University, vol.~3,
  Australian National University, Centre for Mathematical Analysis, Canberra,
  1983.

\end{thebibliography}
\bibliographystyle{amsplain}
\end{document}